\documentclass[12pt]{amsart}
\usepackage{amssymb}
\usepackage{amsxtra}
\usepackage[mathscr]{eucal}
\usepackage{graphics}
\usepackage{epsfig}

\newtheorem{theorem}{Theorem}[section]
\newtheorem{lemma}[theorem]{Lemma}
\newtheorem{corollary}[theorem]{Corollary}

\theoremstyle{definition}

\newtheorem{example}[theorem]{Example}

\def\Per{\operatorname{Per}}
\def\trace{\operatorname{trace}}        
\def\ftrace{\operatorname{T}}   
\def\norm{\operatorname{N}}     
\begin{document}
\title{Integer sequences counting periodic points}
\author{Graham Everest}\email{g.everest@uea.ac.uk}
\author{Yash Puri}\email{yash\underline{\ }puri@hotmail.com}
\author{Tom Ward}\email{t.ward@uea.ac.uk}
\address{School of Mathematics, University of East Anglia,
Norwich NR4 7TJ, UK.} \dedicatory{\today}
\thanks{The second author acknowledges the support of
EPSRC postgraduate award 96001638}

\maketitle

\section{Introduction}

An existing dialogue between number theory and dynamical systems
is advanced.
A combinatorial device gives necessary and
sufficient conditions for a sequence of non-negative integers to
count the periodic points in a dynamical system. This is applied
to study linear recurrence sequences which count periodic points.
Instances where the $p$-parts of an integer
sequence themselves count periodic points are studied.
The Mersenne sequence
provides one example, and the denominators of the Bernoulli numbers
provide another.
The methods
give a dynamical interpretation of many classical congruences such
as Euler-Fermat for matrices, and suggest the same for the
classical Kummer congruences satisfied by the
Bernoulli numbers.

Let $M_n=2^n-1, n\ge1$ denote the $n$-th term of the Mersenne
sequence $\left(M_n\right)$. This sequence is of interest in
number theory because it is expected to contain infinitely many
prime terms, and in dynamics because it counts the periodic points
in the simplest expanding dynamical system. Let $T:{\mathbb
S}^1\rightarrow {\mathbb S}^1$ be the squaring map $T(z)=z^2$, and
let Per$_n(T)$ denote the set of points of period $n$ under $T$,
that is the set of solutions of the equation $T^n(z)=z$. Then it
is easy to check that $\vert\Per_n(T)\vert=M_n$.

Other classical sequences arise in a similar way. Let $L_n$ denote
the $n$-th term of the Lucas sequence $1,3,4,7\ldots$, and let $X$
denote the set of all doubly-infinite strings of $0$'s and $1$'s
in which every zero is followed by a $1$, and $T:X\rightarrow X$
the left shift defined by $(Tx)_n=x_{n+1}$. Then
$\vert\Per_n(T)\vert=L_n$.

The Lehmer-Pierce sequences (generalising the Mersenne sequence;
see \cite{MR2000e:11087}) also arise in counting periodic points.
Let $f(x)$ denote a monic, integral polynomial with degree $d\ge
1$ and roots $\alpha_1,\ldots ,\alpha_d$. Define
$$\Delta_n(f)=\prod_i\vert \alpha_i^n-1\vert,
$$
which is non-zero for $n\ge1$ if no $\alpha_i$ is a root of unity.
When $f(x)=x-2$, we obtain $\Delta_n(f)=M_n$. Sequences of the
form $\left(\Delta_n(f)\right)$ were studied by Pierce and Lehmer
with a view to understanding the special form of their factors, in
the hope of using them to produce large primes. In dynamics they
arise as sequences of periodic points for toral endomorphisms. Let
$X=\mathbb T^d$ denote the $d$-dimensional additive torus. Then
the companion matrix $A_f$ of $f$ acts on $X$ by multiplication
mod $1$, $T(x)=A_fx$ mod $1$. It requires a little thought to
check that $\vert\Per_n(T)\vert=\Delta_n(f)$ under the same {\sl
ergodicity} condition that no $\alpha_i$ is a root of unity (see
\cite{MR2000e:11087}). Notice that the Lehmer-Pierce sequences are
the absolute values of integer sequences which could have mixed
signs.

Our final examples illuminate the same issue of signed sequences
whose absolute value counts periodic points. The Jacobsthal-Lucas
sequence
$R_n=\vert(-2)^n-1\vert$ counts points of period $n$ for the map
$z\mapsto z^{-2}$ on ${\mathbb S}^1$. The sequence $S_n=\vert
2^n+(-3)^n\vert$ counts periodic points in a certain continuous
automorphism of a $1$-dimensional solenoid, see \cite{MR99b:11089}
or \cite{MR90a:28031}.

Following \cite{puri}, call a sequence $u_n$ of non-negative
integers {\sl realisable} if there is a set $X$ and a map
$T:X\rightarrow X$ such that $u_n=\vert\Per_n(T)\vert$. The
examples above were of sequences that are realisable by continuous
maps of compact spaces; it turns out that any realisable sequence
is in fact realisable by such a map.

It is natural to ask what is required of a sequence in order that
it be realisable. For example, could the Fibonacci sequence, the
more illustrious cousin of the Lucas sequence, be realised in this
way? The answer is no, and a simple proof will follow in the next
section. In fact a sequence of non-negative integers satisfying
the Fibonacci recurrence is realisable if and only if it is a
non-negative integer multiple of the Lucas sequence (see
\cite{puri} and \cite{puri-ward}).

The statements of the main theorems now follow. For the first,
note that if $\left(u_n\right)$ is any sequence of integers, then
it is reasonable to ask if the sequence $\left(\vert
u_n\vert\right)$ of absolute values is realisable. For example,
the sequence $1,-3,4,-7,\ldots$ is a signed linear recurrence
sequence whose absolute values are realisable.

The first theorem gives a generalisation of the observation about
realisable sequences which satisfy a linear recurrence relation.
The definitions are standard but they will be recalled in the next
section. Recall that the $\mathbb C$-space of all solutions of a
binary recurrence relation has dimension 2. The {\sl realisable
subspace} is the subspace generated by the realisable solutions.
For the Fibonacci recurrence, the realisable subspace has
dimension 1 and is spanned by the Lucas sequence.

\begin{theorem}\label{binaryrec}
Let $u_n$ denote the $n$-th term of an integer sequence which
satisfies a non-degenerate binary recurrence relation. Let
$\Delta$ denote the discriminant of the characteristic polynomial
associated to the recurrence relation. Then the realisable
subspace has\begin{enumerate}
\item dimension $0$ if $\Delta <0$,
\item dimension $1$ if $\Delta=0$ or $\Delta > 0$
and non-square, and
\item dimension $2$ if $\Delta >0$ is a square.
\end{enumerate}
\end{theorem}

Theorem \ref{binaryrec} surely has a generalisation to higher
degree which characterises the realisable subspace in terms of the
factorisation of the characteristic polynomial $f$. The second
theorem is a partial result in that direction, giving a
restriction on the dimension of the realisable subspace under the
assumption that the characteristic polynomial
has a dominant root.

\begin{theorem}\label{rankrestriction}Let $f$ denote the
characteristic polynomial of a non-degenerate linear recurrence
sequence with integer coefficients. If $f$ is separable and
has $l$ irreducible
factors
and a dominant root then the dimension of the
realisable subspace is
$\leq l$.
\end{theorem}

It is not clear if there is an exact result, but the
deep result of Kim, Ormes and Roush \cite{kor} on
the Spectral Conjecture of Boyle and Handelman \cite{bh}
gives a checkable criterion for a given linear recurrence sequence
to be realised {\sl by an irreducible subshift of finite type}.

\begin{example}\label{tribonacci}
Consider the sequences which satisfy the Tribonacci relation
\begin{equation}\label{trib}
u_{n+3}=u_{n+2}+u_{n+1}+u_n.
\end{equation}
\end{example}
The sequence $1,3,7,11,21,\dots$ satisfies (\ref{trib}) and is
realisable. This is the sequence of traces $\ftrace(A^n)$, where
$A$ is the companion matrix to $f(x)=x^3-x^2-x-1$,
$$A=\left(
\begin{matrix}0&1&0\\0&0&1\\1&1&1\end{matrix}\right).
$$
For an explanation of this remark, turn to the proof of
Corollary
\ref{matrixeulerfermat}. Theorem \ref{rankrestriction} says that
any realisable sequence which satisfies (\ref{trib}) is a multiple
of this one.

The third theorem consists of a pair of examples. Given a
sequence $u_n$ and a prime $p$, write $[u_n]_p$ for the
$p$-part of $u_n$. We say a sequence is {\sl locally
realisable at p} if the sequence $[u_n]_p$ is itself
realisable. We say the sequence is {\sl everywhere locally
realisable} if it is locally realisable at $p$ for all primes
$p$. If a sequence is everywhere locally realisable,
then for each $n\ge1$, $[u_n]_p=1$ for all
but finitely many $p$, and
it is realisable by Corollary \ref{sumandtimes}. We will
sometimes use the term globally realisable for a sequence
when we wish to emphasize the distinction with local realisability.
Consider the Bernoulli numbers,
which are defined by the formula
$$\frac{t}{e^t-1}=\sum_{n=0}^{\infty}B_n\frac{t^n}{n!}.$$
Then $B_n\in \mathbb Q$ for all $n$, and $B_n=0$ for all odd
$n>1$.

\begin{theorem}\label{bernoulliandmersenne}
\begin{enumerate}
\item Any Lehmer--Pierce sequence
is everywhere locally, and hence globally,
realisable.
\item Let $b_n$ denote the denominator of $B_{2n}$ for $n\ge 1$. Then
$\left(b_n\right)$ is everywhere locally, and hence globally,
realisable.
\end{enumerate}
\end{theorem}

The maps in Theorem \ref{bernoulliandmersenne} are
endomorphisms of groups.
Theorem \ref{bernoulliandmersenne} and Lemma \ref{puriwardlemma} suggest a
dynamical interpretation of composite versions of the classical
Kummer congruences; see section \ref{bernoullisect} below.

\section{Combinatorial dynamics}

We begin with a simple remark that shows the Fibonacci sequence is
not realisable. No map can have 1 fixed point and exactly
2 points of
period 3, as any point of least period 3 must have
an orbit of length 3 comprising points all of period 3.
More generally, for any prime $p$, the
number of non-fixed points of period $p$ must be divisible by $p$
because their orbits occur in cycles of length $p$. Using a
generalisation of this kind of reasoning, the following
characterisation emerges.

\begin{lemma}\label{puriwardlemma} Let $u=\left(u_n\right)$ be
a sequence of non-negative integers, and let $u*\mu$ denote the
Dirichlet convolution of $u$ with the M\"obius function $\mu$.
Then $u$ is realisable if and only if $(u*\mu)(n)\equiv 0$ mod $n$
and $(u*\mu)(n)\ge 0$ for all $n\ge 1$.
\end{lemma}

To see why this holds, notice that the set of points of
period $n$ is the {\sl disjoint union} of the set of
points of least period $d$ for $d$ running through the divisors of
$n$, and the number of points with least period $d$
is a mutliple of $d$.
The Dirichlet convolution is the usual definition from analytic
number theory: $\mu(1)=1$, $\mu(n)=0$ unless $n$ is square-free
and $\mu(n)=(-1)^k$ if $n$ is the product of $k$ distinct primes,
and the Dirichlet convolution of an arithmetical function $g$ with
$\mu$ is given by
$$(g*\mu)(n)=\sum_{d\vert n}\mu(d)g(n/d).
$$
Finally, the result is obtained using the M\"obius inversion
formula.
For brevity, write $u^*_n=(u*\mu)(n)$ for $n\ge1$.

\begin{corollary}\label{sumandtimes}
The sum and product of two realisable sequences are both
realisable.
\end{corollary}

This may be seen either using elementary properties of the
Dirichlet convolution or using the realising maps: if $u$ and $v$
are realisable, then the Cartesian product of the realizing maps
realises $\left(u_nv_n\right)$, while the disjoint union realises
$\left(u_n+v_n\right)$.

Notice that if $n=p^r$, for a prime $p$ and $r>0$ an integer,
Lemma \ref{puriwardlemma} requires that
$$u_{p^r}\equiv u_{p^{r-1}} \mod p^r
$$
for any realisable sequence $u$.

\begin{corollary}\label{eulerfermat}
Let $a$ denote a positive integer and let $p$ and $r$ be as above.
Then
$$a^{p^r}\equiv a^{p^{r-1}} \mod p^r.
$$
\end{corollary}

\begin{proof}
This is the statement of the Euler-Fermat Theorem, which may be
seen because the sequence $u_n=a^n$ is realisable. For example,
the left shift $T$ on $\{0,1,\dots,a-1\}^{\mathbb Z}$ has
$\vert\Per_n(T)\vert=a^n$.
\end{proof}

This kind of observation --- that periodic
points in full shifts give simple proofs of many elementary
congruences --- is folklore; indeed the paper \cite{!!!} gives a
rather complicated proof of Euler--Fermat using a dynamical
system.

Lemma \ref{puriwardlemma} does more with no additional effort. The
following is a generalisation of the Euler-Fermat Theorem for
integral matrices which will be used in the proof of Theorem
\ref{binaryrec}.

\begin{corollary}\label{matrixeulerfermat}
Let $A$ denote a square matrix with integer entries and let $p$
and $r$ be as above. Then
$$
\trace(A^{p^r})\equiv\trace(A^{p^{r-1}}) \mod p^r.
$$
\end{corollary}
\begin{proof}
It is sufficient to assume $A$ has non-negative entries, since any
matrix has such a representative mod $p^r$. For non-negative
entries, $\left(\trace(A^n)\right)$ is realisable: Let $G_A$ be
the labeled graph with adjacency matrix $A$ and $T_A$ the
edge-shift on the set of labels of infinite paths on $G_A$. Then
the number of points of period $n$ for this system is
$\trace(A^n)$ (see \cite{LM} for the details).
\end{proof}

We now state the consequences of Lemma \ref{puriwardlemma} in
their most general form for matrix traces.

\begin{corollary}\label{matrixcong}
Let $A$ denote a square matrix with integer entries and let $A_n$
denote the sequence $\trace(A^n)$. Then for all $n\ge 1$
$$A_n^* \equiv 0 \mod n.
$$
\end{corollary}

Before the proof of Theorem \ref{binaryrec}, we begin with some
notation (for a lively account of the general properties of linear
recurrence sequences, see \cite{MR92k:11011}). Suppose we are
given a binary recurrence sequence $u=(u_n)$. This means that
$u_1$ and $u_2$ are given as initial values with subsequent terms
defined by a recurrence relation
\begin{equation}\label{recrel}
u_{n+2}=Bu_{n+1}-Cu_n.
\end{equation}
The polynomial $f(x)=x^2-Bx+C$ is the {\sl characteristic
polynomial} of the recurrence relation. We will write
$$A_f=\left(\begin{matrix}0&1\\ -C&B\end{matrix}\right)
$$
for the companion matrix of $f$. The zeros $\alpha_1$ and
$\alpha_2$ of $f$, are the {\sl characteristic roots} of the
recurrence relation. The assumption on non-degeneracy means that
$\alpha_1/\alpha_2$ is not a root of unity. The {\sl discriminant}
of the recurrence relation is $\Delta = B^2-4C$. Of course, if
$\Delta =0$ then the roots of $f$ coincide, if $\Delta <0$ the
roots are non-real and distinct, if $\Delta>0$ is a square then
the roots are rational and in the other case, the roots are real
and distinct but irrational.

The general solution of the recurrence relation in these cases is
as follows:

$\Delta=0$: $u_n=(\gamma_1+\gamma_2n)\alpha_1^n$ (here
$\alpha_1=\alpha_2$).

$\Delta \neq 0$: $u_n=\gamma_1\alpha_1^n+\gamma_2\alpha_2^n$.

\begin{proof} (of Theorem \ref{binaryrec})

Assume first that $\Delta=0$, and let $p$ denote any prime which
does not divide $\gamma_2$. Then the congruence in Corollary
\ref{matrixeulerfermat} is plainly violated at $n=p$ unless
$\gamma_2=0$. In that case, $\vert\gamma_1\alpha_1^n\vert$ is
realisable and the space this generates is 1-dimensional.

If $\Delta>0$ is a square, then the roots are rationals and
plainly, must be integers. We claim that for any integers
$\gamma_1$ and $\gamma_2$, the sequence $\vert\gamma_1\alpha_1^n
+\gamma_2\alpha_2^n\vert$ is realisable. In fact (up to
multiplying and adding full shifts) this sequence counts the periodic
points for an automorphism on a one-dimensional solenoid, see
\cite{MR2000e:11087} or \cite{MR90a:28031}.

The two cases where $\Delta\neq 0$ is not a square are similar.
Write $\alpha=e+f\sqrt {\Delta}$ for one of the roots of $f$ and
let $K=\mathbb Q(\alpha)$ denote the quadratic number field
generated by $\alpha$. Write $\ftrace_{K\vert \mathbb
Q}:K\rightarrow \mathbb Q$ for the usual field trace. The general
integral solution to the recurrence is $u_n=\ftrace_{K\vert
\mathbb Q}((a+b\sqrt {\Delta})\alpha^n)$, where $a$ and $b$ are
both integers or both half-odd integers. Write
$v_n=\ftrace_{K\vert \mathbb Q}(a\alpha^n)$ and
$w_n=\ftrace_{K\vert \mathbb Q}(b\sqrt {\Delta}\alpha^n)$. Now
$v_n=a\trace(A_f^n)$, where $A_f$ denotes the companion matrix of
$f$. Hence it satisfies $v_p\equiv v_1$ mod $p$ for all primes $p$
by Corollary \ref{matrixeulerfermat}.

Let $p$ denote any inert prime for $K$. The residue field is
isomorphic to the field $\mathbb F_{p^2}$. Moreover, the
non-trivial field isomorphism restricts to the Frobenius at the
finite field level. Reducing mod $p$ gives the congruence
$$
w_p=\ftrace\left(\sqrt{\Delta}\alpha\right)\equiv
\sqrt {\Delta}\alpha^p-\sqrt {\Delta}\alpha\mod p.
$$
Thus, $w_p\equiv -w_1$ mod $p$ for all inert primes $p$. On the
other hand, $v_p\equiv v_1$ mod $p$ for all inert primes $p$.

If $\vert u_n\vert$ is realisable then $\vert u_p \vert \equiv
\vert u_1\vert$ mod $p$ by Corollary \ref{eulerfermat}. If
$u_p\equiv -u_1$ mod $p$ for infinitely many primes $p$ then
$v_p+w_p\equiv v_1-w_1\equiv -v_1-w_1$ mod $p$. We deduce that
$p\vert v_1$ for infinitely primes and hence $v_1=2ae=0$. We
cannot have $e=0$ by the non-degeneracy, so $a=0$. If $u_p\equiv
u_1$ mod $p$ then, by a similar argument, we deduce that $bf=0$.
We cannot have $f=0$ again, by the non-degeneracy so $b=0$. This
proves that when $\Delta \neq 0$ is not a square, the realisable
subspace must have rank less than 2.

Suppose firstly that $\Delta>0$. We will prove that the rank is
precisely 1. In this case, there is a dominant root. If this root
is positive then all the terms of $u_n$ are positive. If the
dominant term is negative then the sequence of absolute values
agrees with the sequence obtained by replacing $\alpha$ by
$-\alpha$ and the dominant root is now positive. In the recurrence
relation (\ref{recrel}) $C=\norm_{K\vert \mathbb Q}(\alpha)$, the
field norm, and $B= \ftrace_{K\vert \mathbb Q}(\alpha)$. We are
assuming $B>0$. If $C<0$ then the sequence $u_n=\ftrace (A_f^n)$
is realisable because the matrix $A_f$ has non-negative entries.
If $C>0$ then we may conjugate $A_f$ to such a matrix (this leaves
the sequence of traces invariant). To see this, let $E$ denote the
matrix
$$E=\left(\begin{matrix}1&0\\k&1\end{matrix}\right).
$$
Then
$$E^{-1}A_fE=\left(\begin{matrix}k&1\\Bk-k^2-C&B-k\end{matrix}\right).
$$
If $B$ is even, take $k=B/2$. Then the lower entries in
$E^{-1}A_fE$ are $(B^2-4C)/4=\Delta/4>0$ and $B/2>0$. If $B$ is
odd, take $k=(B+1)/2$. Then the lower entries are $(B^2-1-4C)/4=
(\Delta-1)/4\geq 0$ and $(B-1)/2\ge0$. In both cases we have
conjugated $A_f$ to a matrix with non-negative entries. Since we
know that the sequence of traces of a matrix with non-negative
entries is realisable, we have completed this part of the proof.

Finally, we must show that when $\Delta<0$, both sequences $v_n$
and $w_n$ are not realisable in absolute value. Assume $a\neq 0$,
and then note that $v_1=2ae\neq0$ by the non-degeneracy
assumption. For all primes $p$ we have $v_p\equiv v_1$ by the
remark above. Since the roots $\alpha_1$ and $\alpha_2$ are
complex conjugates, $\vert\alpha_1\vert=\vert\alpha_2\vert$. Let
$\beta=\frac{1}{2\pi}\arg(\alpha_1/\alpha_2)$; $\beta$ is
irrational by the non-degeneracy assumption. The sequence of
fractional parts of $p\beta$, with $p$ running through the primes,
is dense in $(0,1)$ (this was proved by Vinogradov \cite{vino};
see \cite{vaughan} for a modern treatment). It follows that there
are infinitely many primes $p$ for which $v_pv_1<0$. Therefore, if
$\vert v_n\vert$ is realisable then it satisfies $v_p\equiv v_1$
mod $p$ and $-v_p \equiv v_1$ mod $p$ for infinitely many primes.
We deduce that $v_1=0$ which is a contradiction. With $w_n$ we may
argue in a similar way to obtain a contradiction to $w_1 \neq 0$.
If $\vert w_n\vert$ is realisable then Lemma \ref{puriwardlemma}
says $\vert w_{p^2}\vert \equiv \vert w_p \vert \equiv \vert w_1
\vert$ for all primes $p$. Arguing as before, $w_{p^2}\equiv w_1$
for both split and inert primes. However, the sequence
$\{p^2\beta\}$, $p$ running over the primes, is dense in $(0,1)$.
(Again, this is due to Vinogradov in \cite{vino} or see
\cite{ghosh} for a modern treatment. The general case of
$\{F(p)\}$, where $F$ is a polynomial can be found in
\cite{harman}.) We deduce that $w_{p^2}w_1<0$ for infinitely many
primes. This means $w_{p^2}\equiv w_1$ mod $p$ and $w_{p^2}\equiv
- w_1$ mod $p$ infinitely often. This forces $w_1=0$ - a
contradiction.
\end{proof}

\section{Proof of Theorem \ref{rankrestriction}}

The proof of Theorem \ref{rankrestriction} uses the methods
introduced in the proof of Theorem \ref{binaryrec}.

\begin{proof}
Let $d$ denote the degree of $f$. In the first place
we assume $l=1$, thus $f$ is irreducible.
The irreducibility of $f$ implies
that the rational solutions of the recurrence are given by
$u_n=\ftrace_{K\vert \mathbb Q}(\gamma \alpha^n)$, where
$K=\mathbb Q(\alpha)$, and $\gamma \in K$. We write $\gamma_i,
\alpha_i, i=1,\dots ,d$ for the algebraic conjugates of $\gamma$
and $\alpha$. The dominant root hypothesis says, after
re-labelling, $\vert\alpha_1\vert> \vert \alpha_i\vert$ for
$i=2,\dots ,d$. We will show that if $u_n$ is realisable then
$\gamma \in \mathbb Q$.

Let $p$ denote any inert prime. If $p$ is sufficiently large, the
dominant root hypothesis guarantees that $u_p,\dots,u_{p^{d}}$
will all have the same sign. Using Lemma \ref{puriwardlemma}
several times, we deduce that
$$
u_p\equiv u_{p^2}\equiv \dots \equiv
u_{p^{d}} \equiv \pm u_1 \mod p.
$$
Therefore $u_p+\dots +u_{p^{d}}\equiv \pm du_1$ mod $p$,
the sign depending upon the sign of $u_1$. However,
$$u_p+\dots +u_{p^{d}} \equiv
\ftrace_{K\vert \mathbb Q}(\gamma)\ftrace_{K\vert \mathbb
Q}(\alpha) \mod p.
$$
We deduce a fundamental congruence
$$\ftrace_{K\vert \mathbb Q}(\gamma)\ftrace_{K\vert \mathbb Q}(\alpha)
\equiv \pm d\ftrace_{K\vert \mathbb Q}(\gamma \alpha) \mod
p.
$$
Since this holds for infinitely many primes $p$, the congruence is
actually an equality,
\begin{equation}\label{fundeq}
\ftrace_{K\vert \mathbb Q}(\gamma)\ftrace_{K\vert \mathbb
Q}(\alpha) = \pm d\ftrace_{K\vert \mathbb Q}(\gamma \alpha).
\end{equation}

The next step comes with the observation that if $u_n$ is
realisable then $u_{rn}$ is realisable for every $r\ge 1$. Thus
equation (\ref{fundeq}) now reads
\begin{equation}\label{genfundeq}
\ftrace_{K\vert \mathbb Q}(\gamma)\ftrace_{K\vert \mathbb
Q}(\alpha^r) = \pm d\ftrace_{K\vert \mathbb Q}(\gamma \alpha^r).
\end{equation}
Dividing equation (\ref{genfundeq}) by
$\alpha_1^r$ and letting $r\rightarrow \infty$ we obtain the
equation
$$\ftrace_{K\vert \mathbb Q}(\gamma)= \pm
d\gamma_1.
$$
This means that one conjugate of $\gamma$ is rational and hence
$\gamma$ is rational.

The end of the proof in the case $l=1$ can be re-worked in
a way that makes it more amenable to generalisation. The trace
is a $\mathbb Q$-linear map on $K$ so its kernel has rank $d-1$.
Thus every element $\gamma$ of $K$ can be written $q+\gamma_0$
where $q\in \mathbb Q$ and $\ftrace_{K\vert \mathbb Q}(\gamma_0)=0$.
Noting that $\ftrace_{K\vert \mathbb Q}(q)=dq$ and cancelling
$d$, this simply means equation (\ref{genfundeq}) can be written
$$u_r=\pm
q\ftrace_{K\vert \mathbb Q}(\alpha^r),
$$
for all $r\ge 1$ confirming that the realisable subspace has
rank $\le 1$.

The general case is similar. Each of the irreducible factors
of $f$ generates a number field $K_j, j=1,\dots ,l$ of degree
$d_j=[K_j:\mathbb Q]$. The solutions
of the recurrence look like
$$u_n=\sum_{j=1}^l\ftrace_{K_j\vert \mathbb Q}(\gamma_j\alpha_j^n),
$$
where each $\gamma_j \in K_j$.
Let $L$ denote the compositum of the $K_j$. Using the inert primes
of $L$ and noting that each is inert in each $K_j$, we deduce an
equation
\begin{equation}\label{gengenfundeq}
\sum_{j=1}^l\frac{d}{d_j}
\ftrace_{K_j\vert \mathbb Q}(\gamma_j)
\ftrace_{K_j\vert \mathbb Q}(\alpha_j)
=\pm d
\sum_{j=1}^l\ftrace_{K_j\vert \mathbb Q}(\gamma_j\alpha_j).
\end{equation}
As before, replace $\alpha_j$ by $\alpha_j^r$, and cancel $d$ so that
$$u_r=\pm
\sum_{j=1}^l\frac{1}{d_j}
\ftrace_{K_j\vert \mathbb Q}(\gamma_j)
\ftrace_{K_j\vert \mathbb Q}(\alpha_j^r)
$$
Each $\gamma_j$ can be written $\gamma_j=q_j+\gamma_{0j}$,
where $\ftrace_{K_j\vert \mathbb Q}(\gamma_{0j})=0$. Noting
that $\ftrace_{K_j\vert \mathbb Q}(q_j)=d_jq_j$
we deduce
that
$$u_r=\pm
\sum_{j=1}^lq_j
\ftrace_{K_j\vert \mathbb Q}(\alpha_j^r)
$$
which proves that the realisable subspace has rank $\le l$.
\end{proof}

\section{Proof of Theorem \ref{bernoulliandmersenne}}\label{bernoullisect}

It is sufficient to construct
local maps $T_p:X_p\rightarrow X_p$ for each prime $p$.
Then Corollary \ref{sumandtimes} guarantees a global
realisation by defining
$$T=\prod _pT_p \mbox{ on } X=\prod _pX_p.
$$
If the maps $T_p$ are group endomorphisms then the map $T$ is
a group endomorphism.
\begin{proof}
As motivation, consider
the Mersenne sequence.
For each prime $p$, let
$\mathbb U_p\subset\mathbb S^1$ denote the group of all
$p$th power roots of unity. Define the local endomorphism
$S_p:x\mapsto x^2$ on $\mathbb U_p$. Then $\vert\Per_n(S_p)
\vert=[2^n-1]_p$ so $S_p$
gives a local realisation of the Mersenne sequence.

An alternative approach is to use the $S$-integer dynamical
systems from \cite{MR99b:11089}: for each prime $p$, define
$T_p$ to be the automorphism dual to $x\mapsto 2x$
on $\mathbb Z_{(p)}$ (the localisation at $p$).
Then by \cite{MR99b:11089},
$$
\vert\Per_n(T_p)\vert=\prod_{q\le\infty;q\neq p}
\vert 2^n-1\vert_q=[2^n-1]_p
$$
by the product formula.
This approach gives a convenient proof of the general case.
We may assume that the polynomial $f$ is irreducible;
let $K=\mathbb Q(\xi)$ for some zero of $f$.
Then for each prime $p$, let $S$ comprise all places of
$K$ except those lying above $p$, and let $T_p$ be the
$S$-integer map dual to $x\mapsto \xi x$ on the ring
of $S$-integers in $K$.
Then by the product formula
$$
\vert\Per_n(T_p)\vert=\Big(\prod_{v\vert p}\vert\xi^n-1\vert_v
\Big)^{-1}=\left[\Delta_n(f)\right]_p
$$
as required.

For the Bernoulli denominators,
define
$X_p =\mathbb F_p=
\mathbb Z/p\mathbb Z$. For $p=2$ define $T_p$ to
be the identity. For $p>2$, let $g_p$ denote an
element of (multiplicative) order $(p-1)/2$.
Define $T_p:X_p\rightarrow
X_p$ to be the endomorphism $T_p(x)=g_px$ mod $p$. Plainly
$\vert$Per$_n(T_p)\vert=p$ if and only if $p-1\vert 2n$;
for all other $n$, $\vert$Per$_n(T_p)\vert=1$. The von Staudt--Clausen
Theorem (\cite{MR81i:10002}, \cite{koblitz}) states
that
$$B_{2n}+\sum\frac{1}{p}\in \mathbb Z,
$$
where the sum ranges over primes $p$ for which $p-1\vert 2n$. Thus
$\vert$Per$_n(T_p)\vert= \max\{1,\vert B_{2n}\vert _p\}$ and
this shows the local realisability of the Bernoulli denominators.
\end{proof}

The following statements
seem plausible upon numerical investigation.
\begin{enumerate}
\item the denominators
$$12,120,252,240,132,32760,\dots$$
of $B_{2n}/2n$ form a sequence that is everywhere
locally realisable;
\item the numerators
$$1,1,1,1,1,691,1,3617,43867,\dots$$
of $B_{2n}/2n$ form a realisable sequence
that is not locally realisable at the
irregular primes $37,59,67,101,103,131,149,157,\dots$.
\item the denominators
$$24,240,504,480,264,65520,24,16320,\dots$$
of $B_{2n}/4n$ form a realisable sequence
that is not locally realisable at the primes
$2,3,5,7,11,13$ but seems to be locally realisable for
large primes.
\end{enumerate}

Taking these remarks together with
$n=p^r$ in Lemma \ref{puriwardlemma}, suggests a dynamical
interpretation of the Kummer congruences. These are
stated now, for a proof see \cite{koblitz}.

\begin{theorem}\label{kummercongruences}
If $p$ denotes a prime and $p-1$ does not divide $n$ then $n\equiv
n'$ mod $(p-1)p^r$ implies
$$(1-p^{n-1})\frac{B_n}{n}\equiv (1-p^{n'-1})\frac{B_{n'}}{n'}
\mod p^{r+1}.
$$
\end{theorem}


\end{document}